\documentclass[12pt,titlepage]{article}

\usepackage{graphicx,color}
\usepackage{amsmath,amssymb,amsthm}
\newtheorem{thm}{Theorem}
\newtheorem{lem}[thm]{Lemma}
\newtheorem{cor}[thm]{Corollary}

\newtheorem{prop}[thm]{Proposition}

\newtheorem{defi}{Definition}

\newtheorem{remark}{Remark}

\def\beq{\begin{equation}}\def\eeq{\end{equation}}
\def\beqn{\begin{eqnarray}}\def\eeqn{\end{eqnarray}}

\def\qed{\ifhmode\unskip\nobreak\fi\quad\ifmmode\Box\else$\Box$\fi}

\newcommand{\LL}{{\cal L}}
\newcommand{\KG}{{\rm KG}(n,k)}
\newcommand{\SG}{{\rm SG(n,k)}}

\title{Critical subgraphs of Schrijver graphs for the fractional chromatic number}
\author{\hfil Anna Gujgiczer\thanks
{Department of Computer Science and Information Theory,
Faculty of Electrical Engineering and Informatics,
Budapest University of Technology and Economics and
MTA-BME Lend\"ulet Arithmetic Combinatorics Research Group, ELKH, Budapest, Hungary. Research partially supported by the National Research, Development and Innovation Office (NKFIH) grant K--120706 of NKFIH Hungary; email: {\tt gujgicza@cs.bme.hu}
}
\and G\'abor Simonyi\thanks{Alfr\'ed R\'enyi Institute of Mathematics, Budapest, Hungary and
Department of Computer Science and Information Theory,
Faculty of Electrical Engineering and Informatics,
Budapest University of Technology and Economics. Research partially supported by the National Research, Development and Innovation Office (NKFIH) grants K--120706, K--132696 and SNN-135643 of NKFIH Hungary; email: {\tt simonyi@renyi.hu}
}
}

\date{}

\begin{document}

\maketitle

\begin{abstract}
Schrijver graphs are vertex-color-critical subgraphs of Kneser graphs having the same chromatic number. They also share the value of their fractional chromatic number but Schrijver graphs are not critical for that. Here we present an induced subgraph of every Schrijver graph that is vertex-critical with  respect to the fractional chromatic number. These subgraphs turn out to be isomorphic with certain circular complete graphs. We also characterize the critical edges within this subgraph.

\bigskip
\bigskip
\par\noindent
AMS MSC Primary: 05C15, Secondary: 05C60, 05C72
\end{abstract}

\section{Introduction}

Kneser graphs ${\rm KG}(n,k)$ are defined for every pair of positive integers $n,k$ satisfying $n\ge 2k$. Kneser \cite{Kne} observed (using different terminology) that their chromatic number is not more than $n-2k+2$ and conjectured that this upper bound is tight. This was proved by Lov\'asz in his celebrated paper \cite{LLKn} using the Borsuk-Ulam theorem. Soon afterwards Schrijver \cite{Schr} found that a certain induced subgraph ${\rm SG}(n,k)$ of $\KG$, now called Schrijver graph, still has chromatic number $n-2k+2$ and is also vertex-critical for this property, that is, deleting any of its vertices the chromatic number becomes smaller. It is also well-known that the fractional chromatic number of ${\rm KG}(n,k)$ is $\frac{n}{k}$, a consequence of the vertex-transitivity of these graphs and the Erd\H{o}s--Ko--Rado theorem.
Proving a conjecture of Holroyd and Johnson \cite{Hol} Talbot \cite{Tal} gave the exact value of the independence number of Schrijver graphs that easily implies, as already observed in \cite{ST}, that their fractional chromatic number is also $\frac{n}{k}$. Most Schrijver graphs are not vertex-critical for this property (the only exceptions are the trivial cases when $k=1$, $n=2k$, or $n=2k+1$, cf. Corollary~\ref{cor:Brooks} in Section~\ref{sect:proof}) and this suggests the problem of finding critical subgraphs of Scrijver graphs for the fractional chromatic number. In this paper we present such a subgraph for all values of $n$ and $k$ with $n\ge 2k$. These subgraphs, that turn out to be isomorphic to the circular (also called rational) complete graphs $K_{n'/k'}$ for $n'=\frac{n}{\gcd(n,k)}, k'=\frac{k}{\gcd(n,k)}$, (cf. Corollary~\ref{cor:qknk} in Section~\ref{sect:proof}),
are vertex-transitive, so deleting any of their vertices the value of the fractional chromatic number drops to the same smaller value. We also locate the edges of these special subgraphs that are critical for the fractional chromatic number and show that their deletion already results in the same decrease of the fractional chromatic number as the deletion of a vertex.

The paper is organized as follows. In the next section we give the necessary definitions to define the above mentioned vertex-critical subgraph. In Section~\ref{sect:proof} we prove the vertex-criticality of these graphs and also show the above mentioned relation to circular complete graphs. The last section is devoted to characterizing the critical edges of our subgraphs for the fractional chromatic number.

\section{Well-spread subsets and the subgraph $Q(n,k)$}

\begin{defi}\label{KGSG}
For positive integers $n\ge 2k$ the Kneser graph ${\rm KG}(n,k)$ is defined on the vertex set that consists of the ${n\choose k}$ $k$-element subsets of $[n]=\{1,\dots n\}$ with two such subsets forming an edge if and only if they are disjoint. A $k$-subset $X$ of $[n]$ is called $r$-separated if for any two of its elements $x,y$ we have
$r\le |x-y|\le n-r$. The Schrijver graph ${\rm SG}(n,k)$ is the subgraph of ${\rm KG}(n,k)$ induced by vertices representing $2$-separated sets.
\end{defi}

\medskip
\par\noindent
Notice that arranging the elements of the basic set $[n]$ around a cycle, the $r$-separated sets are exactly those any two elements of which have at least $r-1$ elements on both of the two arcs between them on this cycle.

\medskip
\par\noindent
The following theorem is a condensed version of the well-known results in \cite{LLKn, Schr}.

\begin{thm}\label{cond} {\rm (Lov\'asz--Kneser and Schrijver theorem \cite{LLKn, Schr})}
For every $n\ge 2k$ we have $$\chi({\rm SG}(n,k))=\chi({\rm KG}(n,k))=n-2k+2.$$
Moreover, ${\rm SG}(n,k)$ is vertex-color-critical, i.e.,
$$\forall X\in V({\rm SG}(n,k)):\ \chi({\rm SG}(n,k)\setminus \{X\})=n-2k+1.$$
\end{thm}

\medskip
\par\noindent
The graphs ${\rm KG}(n,k)$ and ${\rm SG}(n,k)$ are widely investigated, cf. e.g. \cite{BBraun1, BBraun2, BjLo, PAC, KSq, KScrit1, KScrit2, Meunier, STSzE80} to mention just a few more of the results related to them.

\medskip
\par\noindent
The fractional chromatic number $\chi_f(G)$ of a graph $G$ is the minimum total non-negative weight one can put on the independent sets of $G$ such that for each vertex the independent sets containing it get at least weight $1$ altogether. It is well-known that, denoting the independence number of graph $G$ by $\alpha(G)$, one always has $$\chi_f(G)\ge\frac{|V(G)|}{\alpha(G)}$$
and equality holds whenever the graph is vertex-transitive, see e.g. \cite{SchU} for this and other basic facts about the fractional chromatic number.

\medskip
\par\noindent
The independence number of Kneser graphs is given by the famous Erd\H{o}s--Ko--Rado theorem.

\medskip
\par\noindent
\begin{thm}\label{EKR} {\rm (Erd\H{o}s--Ko--Rado \cite{EKR})}
$$\alpha({\rm KG}(n,k))={{n-1}\choose {k-1}}.$$
Moreover, for $n>2k$ the only independent sets of this size are the ones whose vertices represent $k$-element subsets that all contain a fixed element $i\in [n]$.
\end{thm}

\medskip
\par\noindent
\begin{cor}\label{cor:frKG}{\rm (cf. e.g. \cite{SchU})}
$$\chi_f({\rm KG}(n,k))=\frac{n}{k}.$$
\end{cor}

\medskip
\par\noindent
Holroyd and Johnson \cite{Hol} conjectured that a similar phenomenon to the one expressed by the Erd\H{o}s--Ko--Rado theorem is also true for Schrijver graphs and more generally, for families of $r$-separated sets. Here we state the result only for $r=2$.

\medskip
\par\noindent
\begin{thm}\label{thm:Tal}{\rm (Talbot \cite{Tal})}
$$\alpha({\rm SG}(n,k))={{n-k-1}\choose {k-1}}.$$
Moreover, for $n>2k, n\neq 2k+2$ the only independent sets of this size in ${\rm SG}(n,k)$ are the ones whose vertices represent $k$-element subsets that all contain a fixed element $i\in [n]$.
For $n=2k+2$ other independent sets of this size exist, too.
\end{thm}

\medskip
\par\noindent
Since $|V({\rm SG}(n,k))|=\frac{n}{k}{{n-k-1}\choose {k-1}}$ and obviously $\chi_f({\rm SG}(n,k))\le\chi_f({\rm KG}(n,k))$ the above theorem has the following immediate consequence already noted in \cite{ST}.

\medskip
\par\noindent
\begin{cor}\label{cor:frSG}
$$\chi_f({\rm SG}(n,k))=\frac{n}{k}.$$
\end{cor}

\medskip
\par\noindent
Let $C_n$ denote the cycle on vertex set $[n]$ where the edges are formed by the pairs of vertices $\{i,i+1\}$ for $i\in\{1,\dots,n-1\}$ and $\{1,n\}.$ In particular, the vertices of ${\rm SG}(n,k)$ are exactly the
$k$-size
independent sets of $C_n$. (We will refer to this cycle as the {\em defining cycle} for $\SG$.)

\medskip
\par\noindent
\begin{defi}\label{defi:wsp}
We call a subset $U$ of $C_n$ {\rm well-spread} if for any two sets $A,B\subseteq [n]$ with $|A|=|B|\le n-1$ satisfying that both induce a (connected) path in $C_n$ we have $$\left||A\cap U|-|B\cap U|\right|\le 1.$$ The induced subgraph of ${\rm SG}(n,k)$ on all well-spread $k$-subsets will be denoted by $Q(n,k)$.
\end{defi}

\medskip
\par\noindent
Note that a set $U\subseteq V(C_n)$ need not be separated (that is $r$-separated for some $r\ge 2$) for being well-spread. Moreover, it follows from the definition that $U$ is well-spread if and only if $\overline{U}:=V(C_n)\setminus U$ is well-spread. Since at least one of $U$ and $\overline{U}$ has size at most $n/2$, one of them is always a separated set and both can be separated if and only if both has exactly $n/2$ elements. Therefore if $\gcd(n,|U|)=1$ (and $n > 2$) then exactly one of $U$ and $\overline{U}$ is a separated set.

\medskip
\par\noindent
Now we state a basic property of the graphs $Q(n,k)$.

\medskip
\par\noindent
\begin{prop}\label{prop:relpr}
Let $n\ge 2k$ and $\ell\ge 2$ be any positive integer. Then the graphs $Q(n,k)$ and $Q(\ell n,\ell k)$ are isomorphic.
\end{prop}

\proof
Let $U \subseteq V(C_{\ell n})=[\ell n]$ be a well-spread set of size $\ell k$. Consider the $n$-element sets $A_i\subseteq [\ell n], i\in [\ell n]$ defined by
$$A_i:=\{i,i+1,\dots,i+n-1\},$$
where the addition is intended modulo $\ell n$ (and $0$ is represented by $\ell n$), that is the sets $A_i$ are exactly those subsets of $[\ell n]$ that induce a path of length $n$ in $C_{\ell n}$. The number of pairs in the set $$\{(j, A_i)\colon j\in A_i\cap U\},$$ where $j\in [\ell n]$ and $A_i$ is one of the sets just defined
is $\ell kn$ since each $j\in U$ will appear in exactly $n$ distinct $A_i$'s and $|U|=\ell k$. Since there are $\ell n$ distinct $A_i$'s, this means that if any $A_i$ would contain less than $k$ elements of $U$, then some other $A_{i'}$ should contain more than $k$ resulting in two sets $A_i, A_{i'}$ of the same size, both inducing a path of $C_{\ell n}$ whose size of intersection with $U$ differs by at least $2$. This would contradict the well-spread property of $U$, so this is impossible. The situation is similar if any $A_i$ would contain more than $k$ elements of $U$, therefore we have $$\forall i\colon |A_i\cap U|=k.$$

\smallskip
\par\noindent
This means that $i\in U$ if and only if $i+n$ (mod $\ell n$) $\in U$ for every $i \in V(C_{\ell n})$  (otherwise $|A_i\cap U|=|A_{i+1}\cap U|$ would not be satisfied).
Hence, if we have $X \in V(Q(\ell n,\ell k))$, that is $X$ is a well-spread $(\ell k)$-subset of $[\ell n]$ and we
rotate the defining cycle $C_{\ell n}$ exactly $n$ times then we get a vertex $Y \in Q(\ell n, \ell k)$ that is identical with  $X$. So we always get back the same vertex by a $jn$ rotation (where $0 \leq j \leq \ell-1$) and thus the same graph that we get if we identify all the vertices of $C_{\ell n}$ that are $jn$-rotations of each other for some $1 \leq j \leq \ell-1$ and consider well-spread $k$-element subsets. This proves the statement.
\qed

\medskip
\par\noindent
Note that Proposition~\ref{prop:relpr} implies that $Q(n,k)\cong Q\left(\frac{n}{\gcd(n,k)},\frac{k}{\gcd(n,k)}\right)$, therefore when discussing the properties of $Q(n,k)$ we may assume that $\gcd(n,k)=1$.

\medskip
\par\noindent
Now we can already state our result on the vertex-criticality of $Q(n,k)$ for the fractional chromatic number.

\medskip
\par\noindent
\begin{thm}\label{thm:main}
Assume $n\ge 2k$, $\gcd(n,k)=1$ and let $a$ and $b$ be the smallest positive integers for which $ak=bn-1.$
The graph $Q(n,k)\subseteq \SG$ satisfies the following properties.
\begin{itemize}
\item
$\chi_f(Q(n,k))=\frac{n}{k}=\chi_f(\SG).$

\item
$\forall U\in V(Q(n,k))\ \ \chi_f(Q(n,k)\setminus\{U\})=\frac{a}{b}<\frac{n}{k},$
that is $Q(n,k)$ is vertex-critical for the fractional chromatic number.

\item
$Q(n,k)$ contains an induced subgraph isomorphic to $Q(a,b)$.

\end{itemize}

\end{thm}

\medskip
\par\noindent
We are going to prove this theorem in the next section.

\section{Proof of vertex criticality}\label{sect:proof}

\medskip
\par\noindent
For proving Theorem~\ref{thm:main} it would be enough to show that if $\gcd(n,k)=1$ then the $Q(n,k)$ subgraph is isomorphic to the circular complete graph $K_{n/k}$ (that we will define later in Definition~\ref{defi:circomp}). This would suffice because of the following reasoning.
The fractional chromatic number of $K_{n/k}$ is $n/k$ since it is vertex transitive and has $n$ vertices, while its independence number is $k$ (cf. \cite{HN}). It is known that removing any vertex $x$ from it, the remaining graph $K_{n/k} - \{x\}$ is homomorphically equivalent to $K_{a/b}$ where $a$ and $b$ are the unique solution for $nb-ka = 1$, see Lemma 6.6 in \cite{HN} where a retract of $K_{n/k}$ isomorphic to $K_{a/b}$ is shown.
This would already give all the three statements of Theorem~\ref{thm:main}.
We will show this required isomorphism later in this section after all the necessary lemmas will have been provided that will also be needed for the different proof we give (which also has some essential similarities to the one of Lemma 6.6 in \cite{HN} for $K_{n/k}$). Nevertheless, it focuses on the $Q(n,k)$ representation and on this graph being a subgraph of $SG(n,k)$. In Section~\ref{sect:edges} we will modify our proof to obtain a stronger result that characterizes the (fractionally) critical edges of $Q(n,k)$ as well.

\medskip
\par\noindent
Our argument will need the following alternative characterization of well-spread $k$-subsets.

\medskip
\par\noindent
\begin{lem}\label{lem:wspalt}
A subset $U\subseteq V(C_n)$ is well-spread if and only if whenever $A,B\subseteq V(C_n)$ are two sets inducing a minimal path containing the same number $s\le |U|$ elements of $U$ then $$||A|-|B||\le 1.$$
\end{lem}

\medskip
\par\noindent
Note that the minimality of the paths referred to in the statement of Lemma~\ref{lem:wspalt} does not mean that they have {\em minimum} length. On the other hand, their minimality implies that both $A$ and $B$ must start and end with vertices of $C_n$ belonging to $U$.

\proof
Assume to the contrary that $||A|-|B|| \geq 2$ and w.l.o.g. assume that $|A|-2 \geq |B|$. Then, we can modify the subset $A$ by removing its two extremal (that is starting and ending) vertices and $|A|-|B|-2$ more vertices from one end. This way we obtain a path $A'$ for which $|A'| = |B|$ but $||A'\cap U|-|B\cap U|| \geq 2$ which means that $U$ cannot be a well-spread set.

\smallskip
\par\noindent
For the other direction suppose that $U$ is not well-spread. Then there exist $A$ and $B$ paths in $C_n$ for which $|A| = |B|$ but $||A\cap U|-|B\cap U|| \geq 2$. W.l.o.g. assume, that $|A\cap U| \geq |B\cap U|+2$. We may assume that $A$ induces a minimal path in $C_n$ containing $s:=|A\cap U|$ elements of $U$ because if not, then we can make both $A$ and $B$ shorter so that $|A\cap U|$ does not change while $|B\cap U|$ may only become smaller, so the relations $|A\cap U| \geq |B\cap U|+2$ and $|A|=|B|$ remain valid. Now extend $B$ at both of its ends until it will contain a new element of $U$ at both ends, that is we obtain a $B'$ which induces a minimal path that intersects $U$ in $|B\cap U|+2$ elements. If this number is still less than $s$ than extend $B'$ further (on one end) to make it a minimal path containing exactly $s$ elements of $U$. Since in the first step we extended $B$ at both ends we certainly have $|B'|\ge |A|+2$, so regarding $A$ and $B'$ the condition in the statement cannot hold for our not well-spread set $U$. This completes the proof.
\qed

\medskip
\par\noindent
\begin{lem}\label{lem:euclid}
Well-spread sets are unique up to rotations of the cycle $C_n$. In particular, the graph $Q(n,k)$ is vertex-transitive for any $n$ and $k$ and if $\gcd(n,k)=1$ then $|V(Q(n,k)|=n$.
\end{lem}

\proof
Let $U$ be a well-spread set of size $k$ on the cycle $C_n$. We may (and do) assume $\gcd(n,k)=1$, since the other cases are taken care of by Proposition~\ref{prop:relpr} (and its proof).
\smallskip
\par\noindent
Call $x,y \in U$ neighboring in $U$ if on one of the arcs between them $\nexists z \in U$. By Lemma~\ref{lem:wspalt} if $x,y\in U$ are neighboring, then $x$ and $y$ should be $q_1 := \left\lfloor\frac{n}{k}\right\rfloor$ or $\left\lfloor\frac{n}{k}\right\rfloor + 1$ distance apart (including $x$, but not $y$, that is, they are separated by $q_1-1$ or $q_1$ other elements of the cycle, respectively). If $n=q_1k+r_1$ then we have exactly $r_1$ neighboring pairs of elements in $U$ whose distance is $q_1+1$ and $k-r_1$ neighboring pairs have distance $q_1$.
Remove $q_1-1$ vertices of $C_n$ from the arc between every neighboring pair $x,y \in U$. This way we obtain a shorter cycle $C_{n-(q_1-1)k}$ on which $U$ is still well-spread and (by the observation given after Definition~\ref{defi:wsp}) $\overline{U} := V(C_{n-(q_1-1)k}) \setminus U$ is well-spread as well. On this shorter cycle $U$ is not separated any more, so $\overline{U}$ is a separated set (again by the observation after Definition~\ref{defi:wsp}). Using the notation $n_1:=n-(q_1-1)k=k+r_1$ we have $|\overline{U}|=n_1-|U|=r_1$ and the neighboring elements of $\overline{U}$ are separated by
$q_2:=\left\lfloor\frac{k+r_1}{r_1}\right\rfloor$ or by $q_2-1$ elements of $U$. Now performing the previous removal process with $C_{n_1}$ in the place of $C_n$ and its $r_1$-element subset $\overline{U}$ in place of $U$ is essentially performing a second step of the Euclidean algorithm with $k+r_1$ and $r_1$ (instead of $k$ and $r_1$ but this is not an essential difference). Thus iterating this process we will arrive to a situation where we have a cycle $C_m$ for some smaller $m$ and our current separated set will have only $\gcd(n,k)=1$ element. Clearly a $1$-element set can be placed on the vertex set of $C_m$ in a unique way (up to rotations of the cycle). Noticing that reversing our removal process is completely deterministic, this proves that the original $k$-element well-spread set can also be placed on $C_n$ only in a unique way up to rotations of the cycle. This proves the first statement from which the vertex-transitivity of $Q(n,k)$ immediately follows along with $|V(Q(n,k))|\le n$.
\medskip
\par\noindent
We still have to prove that $\gcd(n,k)=1$ also implies $|V(Q(n,k))|\ge n$. Assume to the contrary that the shortest rotation $f$ of $C_n$ that maps a vertex $X$ of $Q(n,k)$ to itself is shorter than $n$, that is, $f$ is a $t$-fold rotation for some $t<n$. We claim that $\forall i,j\in X$ the arcs $A_i=(i,i+1,\dots,i+t)$ and $A_j=(j,j+1,\dots,j+t)$ contain the same number of elements of $X$. Let $s:=|A_i\cap X|$ and assume for contradiction w.l.o.g. that $s>|A_j\cap X|$. Then add elements $j+t+1,j+t+2,\dots,j+t+h$ to $A_j$ to obtain $A_j'$, where $h$ is the smallest number that makes $|A_j'\cap X|=|A_i\cap X|=s$. Then $A_i$ and $A_j'$ are both minimal arcs containing $s$ elements of $X$. Since $X$ is a separated set, we have $h\ge 2$. Since $|A_i|=|A_j|$ by their definition, this implies $|A_j'|\ge |A_i|+2$ contradicting the well-spreadness of $X$ by Lemma~\ref{lem:wspalt}.

\smallskip
\par\noindent
That means that we must have $|A_i\cap X|=s=|A_j\cap X|$. In other words this means that every $(s-1)$th element of $X$ is $t$ $(=|A_i|-1)$ clockwise rotations away in $C_n$.
As $f$ was chosen to be the shortest rotation mapping $X$ to itself, and since the $n$-fold rotation also maps $X$ to itself, this means that $s-1$ divides $k$, $t$ divides $n$
and the two ratios are equal, that is, $$\frac{k}{s-1}=\frac{n}{t}=:d.$$ But then $d$ is a common divisor of $n$ and $k$, so it must be $1$ contradicting that $f$ is a less-than-$n$-fold rotation.
\qed

\begin{figure}[!htbp]
\centering
\includegraphics[scale=0.9]{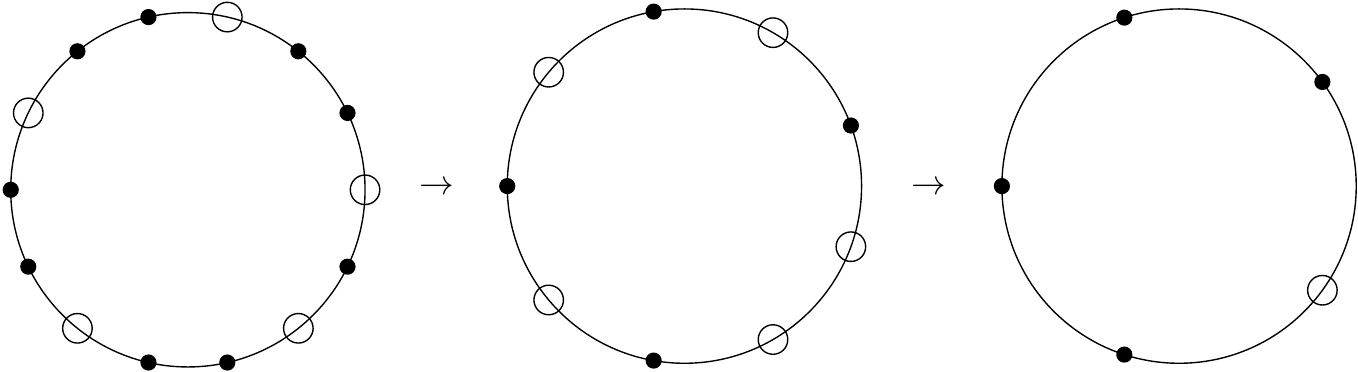}
\caption{The process in the proof of Lemma~\ref{lem:euclid} performed for $n=14, k=5$.}
\label{fig:G145}
\end{figure}

\medskip
\par\noindent
\begin{cor}\label{cor:uniq}
If $\gcd(n,k)=1$ then for every $X,Y\in V(Q(n,k))$ there is a unique rotation of $C_n$ that maps $X$ to $Y$.
\end{cor}

\proof
We have $|V(Q(n,k))| = n$, where the vertices can only be different by some rotation and we have exactly $n$ possible rotations for each vertex.
\qed

\medskip
\par\noindent
\begin{lem}\label{lem:nemfugg}
Let $\gcd(n,k)=1$ and $X,Y\in V(Q(n,k))$ be such that $XY\notin E(Q(n,k))$, that is, $X\cap Y\neq\emptyset$. Let $f\colon V(C_n)\to V(C_n)$ be the unique clockwise rotation moving $X$ to $Y$ and let $i$ be an element of $X\cap Y$. Then the number of elements of $Y$ on the arc of $C_n$ between $i$ and $f(i)$ (moving from $i$ to $f(i)$ in the clockwise direction) is independent of the choice of $i\in X\cap Y$.
\end{lem}

\proof
The argument needed here is almost the same as in the last part of the proof of Lemma~\ref{lem:euclid}.
Let $i, j \in X\cap Y$ and
let $A$ and $B$ be the arcs of $C_n$ between $i$ and $f(i)$ and between $j$ and $f(j)$, respectively ($i,f(i)$ and $j$, $f(j)$ included). We obviously have $|A| = |B|$.
Assume to the contrary that w.l.o.g. $|A \cap Y| + 1 \leq |B \cap Y|$.
Add the minimal number of consecutive vertices to $A$ from $C_n$ in the same (clockwise) direction to get $A'$, such that $|A' \cap Y| = |B \cap Y|$. $Y \in V(SG(n,k))$, therefore it is a 2-separated set, so $|A'| \geq |A| +2 = |B| +2$. Since $A'$ and $B$ are minimal arcs containing the same number of elements of $Y$, this gives a contradiction by Lemma~\ref{lem:wspalt} with the well-spreadness of $Y$.
\qed

\medskip
\par\noindent
\begin{defi}\label{defi:jneighb}
Under the conditions of the previous lemma we call vertex $Y\in V(Q(n,k))$ a right $j$-neighbor of vertex $X\in V(Q(n,k))$ if the number of elements of $Y$ on the arc of $C_n$ between $i$ and $f(i)$ for some $i\in X\cap Y$ (moving from $i$ to $f(i)$ in the clockwise direction) is $j$ (including $f(i)$ but not including $i$)\footnote{Note that the number of elements of $Y$ between $i$ and $f(i)$ is the same as the number of elements of $X$ between $f^{-1}(i)$ and $i$. This observation is behind the terminology we introduce here.}.
\end{defi}

\medskip
\par\noindent
Note that by Lemma~\ref{lem:nemfugg} the previous definition is meaningful as it does not depend on the choice of $i\in X\cap Y$.
Since $X$ is a well-spread set on $C_n$, it has only two different right $j$-neighbors for each $j\in\{1,\dots k-1\}$ and they are exactly one rotation apart from each other (by Lemma~\ref{lem:wspalt}), meaning that they are adjacent since they must be disjoint. Therefore the following is true.

\medskip
\par\noindent
\begin{cor}\label{cor:egyj}
  Let $\gcd(n,k)=1$, $A\subseteq V(Q(n,k))$ be an independent set of $Q(n,k)$ and $X\in A$. Then for each $j\in\{1,\dots k-1\}$ $A$ can contain at most one right $j$-neighbor of $X$. In particular, $$|A|\le k.$$ \hfill$\Box$
\end{cor}

\medskip
\par\noindent
For a complete proof of Theorem~\ref{thm:main} we should investigate further the maximum independent sets of $Q(n,k)$. To this end the following representation of the graph will be helpful that we will call its {\em natural representation}. We arrange the vertices of $Q(n,k)$ around a cycle, which we will call the {\em base cycle}, according to the rotation of the defining cycle of the points of our $n$-set, that is, $X$ is right after $Y$ in the clockwise direction on the base cycle if we get the points of $X$ with one clockwise rotation of the points of $Y$ on the defining cycle. Also, when we talk about an arc {\em starting from} a vertex $X$, we mean, that $X$ is the anticlockwise endpoint of this arc.

\medskip
\par\noindent
\begin{lem}\label{lem:samelook}
In the natural representation of $Q(n,k)$ if we look at any maximum independent set $A \subset V(Q(n,k))$ it forms a well-spread set on the base cycle.
\end{lem}

\proof
Let $A$ be an arbitrary maximum independent set of $Q(n,k)$ and $X_1,X_2\in A$.
Since $A$ is a maximum size independent set both $X_1$ and $X_2$ need to have a right $j$-neighbor in this set for each $j\in\{1,\dots k-1\}$.
Let $Y_1$ be a right $j$-neighbor of $X_1$ and $Y_2$ be a right $j$-neighbor of $X_2$ in $A$ for some $j$. Let $B_1$ be the arc in the cyclic arrangement between $X_1$ and $Y_1$ (containing both $X_1$ and $Y_1$), that is, the number of clockwise rotations needed to get $Y_1$ from $X_1$ is $|B_1|-1$ (both on the defining cycle and on the base cycle), and let $B_2$ be similarly defined for $X_2$ and $Y_2$.
Obviously for $i < j$ the number of clockwise rotations needed for a right $i$-neighbor is less than the number of clockwise rotations needed for a right $j$-neighbor. This means that $|A \cap B_1| = |A \cap B_2|=j+1$. By Lemma~\ref{lem:wspalt} and Definition~\ref{defi:jneighb} we have $||B_1|-|B_2||\leq 1$ (as the distance of right $j$-neighbors on the defining cycle can be one of only two consecutive values) and also by Lemma~\ref{lem:wspalt} this means that $A$ is a well-spread subset of the base cycle.
\qed

\medskip
\par\noindent
The above statements already provide the independence number of $Q(n,k)$ along with an Erd\H{o}s-Ko-Rado type property: the maximum independent sets are formed by $k$-subsets that intersect each other in the same element.

\medskip
\par\noindent
\begin{cor}\label{cor:ekr}
If $\gcd(n,k)=1$, then $$\alpha(Q(n,k))=k$$ and every $k$-element independent set consists of vertices that represent subsets containing a common element.
\end{cor}

\proof
In Corollary~\ref{cor:egyj} we have seen that $\alpha(Q(n,k))\le k$ (we assume $\gcd(n,k)=1$ all along). On the other hand, $k$-element independent sets do exist: consider a vertex $X$ consisting of the elements $i_1,i_2,\dots,i_k$ and then take all its rotated versions on the defining cycle that moves $i_2$ to $i_1$, then $i_3$ to $i_1$, etc., until $i_k$ to $i_1$. This gives a $k$ element independent set each vertex of which contains the element $i_1$. This shows that $\alpha(Q(n,k))\ge k$ and since by Lemma~\ref{lem:samelook} all the maximum independent sets ``look alike'', they all must have the property that each of their vertices contains a common element.
\qed

\medskip
\par\noindent
\begin{lem}\label{lem:frszin}
Let $\gcd(n,k)=1$, $a,b$ be the smallest positive integers for which $ak=bn-1$ and $X\in V(Q(n,k))$. Then $$\chi_f(Q(n,k)\setminus\{X\})\le\frac{a}{b}.$$
\end{lem}

\proof
Arrange the vertices of $Q(n,k)$ along the base cycle according to the natural representation. Let $A_1$ be a maximum independent set in $Q(n,k)$ and let $A_2, \dots, A_a$ be the independent sets we obtain with $1, \dots, a-1$ anticlockwise rotations of $A_1$.
Let $\LL$ be the collection of the $n$ distinct $a$-length arcs of the base cycle.
Since every vertex of $Q(n,k)$ appears in exactly $a$ of the $a$-length arcs of the base cycle and $|A_1|=k$, the set $$\{(v,L)\colon v\in A_1, L\in\LL\}$$ has cardinality $ak$ which by the choice of $a$ and $b$ is equal to $b(n-1) + (b-1)$.
Since by Lemma~\ref{lem:samelook} $A_1$ is well-spread on the base cycle, $|\{(v,L)\colon v\in A_1, L\in\LL\}|=b(n-1)+(b-1)$ means that exactly one of the $n$ arcs of length $a$ in $\LL$ contains $b-1$ elements of $A_1$ and the remaining $n-1$ such arcs contain $b$ elements of $A_1$. (Otherwise we would have $a$-length arcs containing at least $b+1$ and others containing at most $b-1$ elements of $A_1$ contradicting its well-spreadness.) Let $X\in V(Q(n,k))$ be the starting vertex (i.e. anticlockwise endpoint) of the unique $a$-length arc containing only $b-1$ elements of $A_1$. Since $Q(n,k)$ is vertex transitive, we may assume that we delete this vertex $X$. Remove $X$ and put weight $\frac{1}{b}$ on each of the independent sets $A_1, \dots, A_a$. Every vertex $Y \in V(Q(n,k)) \setminus \{ X\}$ will be contained in as many of these independent sets as many elements of $A_1$ are in the $a$-length arc starting with $Y$, which is $b$. Thus the total weight on each such $Y$ is $1$, that is, we obtain a fractional coloring this way with total weight $\frac{a}{b}$ proving the statement.
\qed

\medskip
\par\noindent
{\em Proof of Theorem~\ref{thm:main}.}
\smallskip
\par\noindent
The first statement follows from Lemma~\ref{lem:euclid} and Corollary~\ref{cor:ekr}: $Q(n,k)$ is vertex-transitive on $n$ vertices with independence number $k$, so $\chi_f(Q(n,k))=\frac{n}{k}$.

\medskip
\par\noindent
The second statement of Theorem~\ref{thm:main} follows from Lemma~\ref{lem:frszin} and the third statement so all we have to do is to prove the third statement.

\medskip
\par\noindent
Let $H$ be a subgraph of $Q(n,k)$ induced by any $a$ vertices $X_1, \dots, X_a$ consecutive on the base cycle according to the natural representation. Here we consider $X_{i+1}$ to be one anticlockwise rotation away from $X_i$ (for $i \in \{1,\dots,a-1$\}). We will show that $H$ is isomorphic to $Q(a,b)$.

\smallskip
\par\noindent
Since $X_1$ is a well-spread set on the defining cycle $C_n$, we know by the choice of $a$ and $b$ that there is only one $a$-length arc $B'$ in $C_n$ for which $|B' \cap X_1| = b-1$ (the other $a$-length arcs intersect $X_1$ in $b$ elements). Let $B$ be the arc of $C_n$ that we get from $B'$ with one clockwise rotation. Now, for every $X_i$ ($i \in \{1.\dots,a\}$) we have $|X_i \cap B| = b$ since $a<n$ and only the $(n-1)$-fold anticlockwise rotation of $X_1$ (that we would call $X_n$) has only $b-1$ elements in $B$.

\smallskip
\par\noindent
Let the elements of $B$ be labeled $1,2,\dots,a$ as they follow each other on $C_n$ in the clockwise direction. A simple consequence of the choice of the arc $B$ is that in every rotation moving $X_i$ to $X_{i+1}$ the last vertex (clockwise endpoint) of $B$ will belong to $X_{i+1}$ if and only if the first vertex of $B$ belonged to $X_i$ and this holds for every $i = 1, \dots a-1$. (Otherwise we would have an $a$-length arc different of $B'$ that intersects $X_1$ in a different number of elements than $b$.)
So $B$ "behaves" like an $a$-length defining cycle: if we consider its first endpoint labeled $1$ and last endpoint labeled $a$ as neighbors on an $a$-length cycle, then we find that the $b$-element sets $X_i\cap B$ are all well-spread on the so-obtained $C_a$ and for $1\le i,j\le a$ $X_i\cap X_j=\emptyset$ is equivalent to $(X_i\cap B)\cap (X_j\cap B)=\emptyset$. (The latter can be seen by observing that from the point of view of their relative position rotating $X_i$ one anticlockwise rotation is equivalent to moving $B$ one clockwise rotation.)
Let $g$ be the function that assigns the two endpoints of $B$ to two adjacent vertices of the defining cycle $C_a$ as said above and the rest according to their distance from the endpoints.
Let the vertices of $Q(a,b)$ be $Y_i$ defined by $\{g(v) | v \in X_i \cap B\}$ on $C_a$ for each $i \in [a]$.

\smallskip
\par\noindent
Now we have to show, that the function $g(X_i) = Y_i$ ($i=1,\dots,a$) is an isomorphism.
However, this follows from what is already said above that $X_i\cap X_j=\emptyset$ is equivalent to $(X_i\cap B)\cap (X_j\cap B)=\emptyset$
\qed

\begin{figure}[!htbp]
\centering
\includegraphics[scale=0.9]{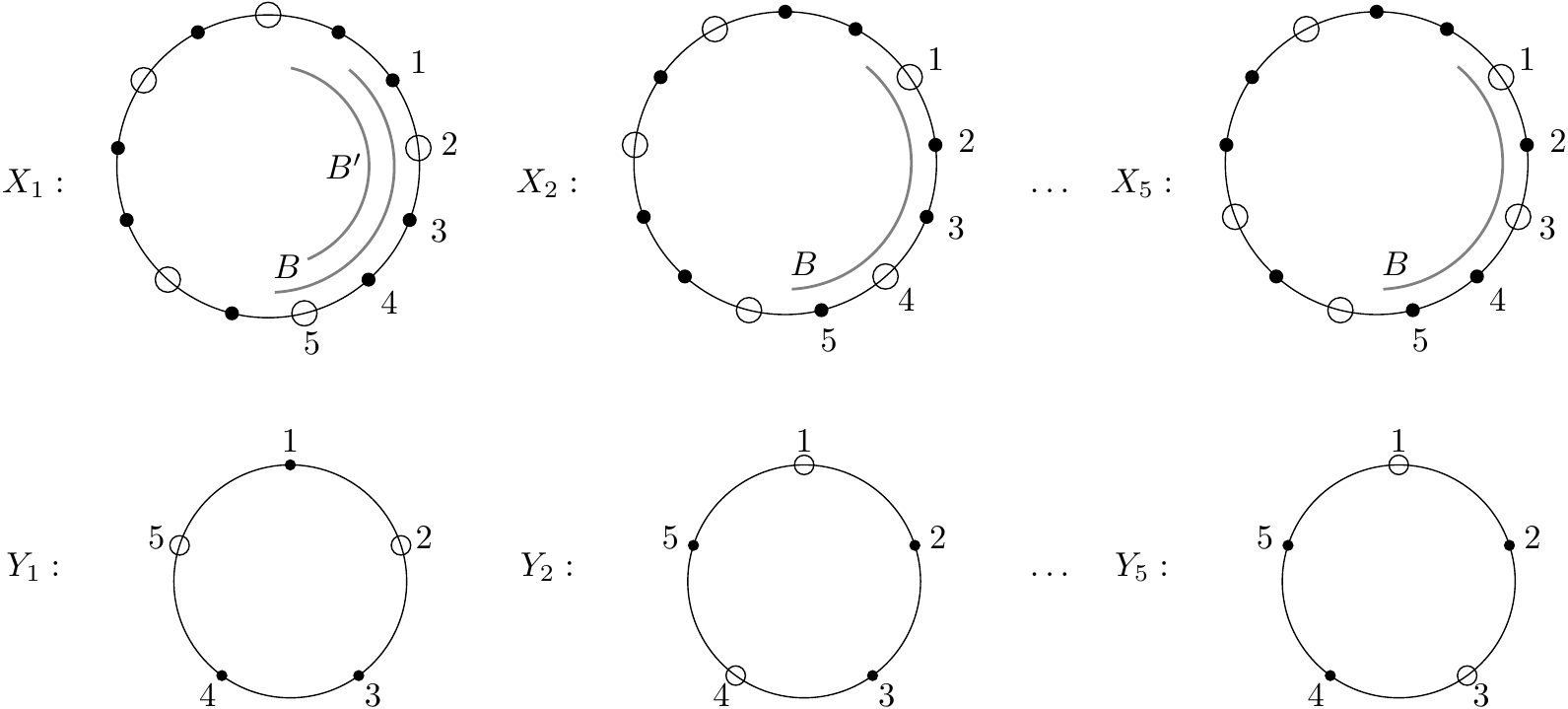}
\caption{Finding the subgraph $Q(5,2)$ in $Q(13,5)$ according to the proof of Theorem~\ref{thm:main}.}
\label{fig:G135}
\end{figure}

\begin{prop}\label{prop:retract}
Assume $n\ge 2k$, $\gcd(n,k)=1$ and let $a$ and $b$ be the smallest positive integers for which $ak=bn-1.$ Then
$\forall U\in V(Q(n,k))\ \ Q(n,k)\setminus\{U\}$
is homomorphically equivalent to $Q(a,b)$.
\end{prop}

\proof
From the previous proof we know that $Q(a,b)$ is a subgraph of $Q(n,k)$ so we have $Q(a,b) \to Q(n,k)$.
By modifying the previous proof we can also easily get the reverse homomorphism. Let arc $B$ and vertex $X_1$ be defined similarly as in the proof of Theorem~\ref{thm:main} above and let $X_2,\dots, X_{n-1}$ be the vertices obtained by $1,2,\dots,n-2$ anticlockwise rotations of $X_1$ on the defining cycle $C_n$. We assume that the deleted vertex $U$ is just the $n^{\rm th}$ vertex of $Q(n,k)$ that is not appearing on this list, that is, it is the vertex we would obtain from $X_1$ by $n-1$ anticlockwise
rotations (or equivalently, with $1$ clockwise rotation). Thus $U$ is the unique vertex for which $|U\cap B|=b-1$, so we still have that for each $i \in [n-1]$ $|X_i \cap B| = b$. This means that whenever some element of $B$ is in $X_i$ for some $i$ then it is also in $X_{i+a}$ (since $|B|=a$).
\smallskip
\par\noindent
Using again the notation $Y_i$ for the vertex of $Q(a,b)$ represented by $X_i\cap B$ if the interval $B$ is ``closed'' to an $a$-length cycle by making its first and last endpoints consecutive elements of $C_a$, the function $$g: X_{i'}\mapsto Y_i\ \ {\rm iff}\ \ i' \equiv i\ ({\rm mod}\ a)$$ is a homomorphism from $Q(n,k)\setminus \{U\}$ to $Q(a,b)$.
\qed

\bigskip

\medskip
\par\noindent
Now we show that ${\rm SG}(n,k)$ itself is critical for the fractional chromatic number only in the cases already mentioned in the Introduction. To this end we first prove the following. (Note that $Q(n,k)$ is a circulant graph, in particular, it is regular.)

\begin{prop}\label{prop:deg}
Let $n>2k$ and assume $\gcd(n,k)=1$. Then the degree of vertices in $Q(n,k)$ is $n-2k+1$.
\end{prop}

\proof
We already know, that a vertex $X \in V(Q(n,k))$ has $2(k-1)$ other vertices not adjacent to it, because  $X$ has exactly two different right $j$-neighbors for each $j\in\{1,\dots k-1\}$. So we have that $X$ is adjacent to exactly $n-2(k-1)-1$ other vertices.
\qed

\begin{cor}\label{cor:Brooks}
We have $Q(n,k)={\rm SG}(n,k)$ if and only if $k=1, n=2k$, or $n=2k+1$. In particular, $\SG$ is vertex-critical for the fractional chromatic number in exactly these cases.
\end{cor}

\proof
We know from Schrijver's theorem, that $\chi(\SG)=n-2k+2$. By Proposition~\ref{prop:deg} this is exactly one more than the (maximum) degree of $Q(n,k)$. Thus (since $\SG$ is connected) by Brooks' theorem, in case $\SG=Q(n,k)$ we must have that $\SG$ is a complete graph or an odd cycle. This happens only in the cases listed in the statement and in those cases we indeed have $Q(n,k)={\rm SG}(n,k)$.
\qed

\medskip
\par\noindent

Now we have all the necessary lemmas to prove that our $Q(n,k)$ graph is isomorphic to the circular complete graph $K_{n/k}$ whenever $\gcd(n,k)=1$.
\medskip
\par\noindent
\begin{defi} \label{defi:circomp}
The circular complete graph $K_{n/k}$ is defined as follows:
$$V(K_{n/k}) = \{0,1, \dots , n-1 \}$$
$$E(K_{n/k}) = \{ \{i,j\} : k \leq |i-j| \leq n-k. \}$$
\end{defi}

\par\noindent
Some important properties of these graphs are that they are vertex-transitive, that $K_{n/k}$ is homomorph equivalent to $K_{n'/k'}$ whenever $\frac{n}{k} = \frac{n'}{k'}$ and that $\chi(K_{n/k}) = \left\lceil\frac{n}{k}\right\rceil$ (for these and further properties, see \cite{HN}). Note that the just stated homomorph equivalence cannot be an isomorphism if $n\neq n'$ since then $|V(K_{n/k})|=n\neq n'=|V(K_{n'/k'})|$. This is a crucial difference between the graphs $K_{n/k}$ and $Q(n,k)$ and shows that the condition $\gcd(n,k)=1$ cannot be dropped in the following statement.

\begin{prop}\label{prop:isotocirc}
$Q(n,k)$ is isomorphic with the circular complete graph $K_{n/k}$ whenever $\gcd(n,k)=1$.
\end{prop}

\proof
As $|V(Q(n,k))| = |V(K_{n/k})| = n$ and in both graphs each vertex has degree $n-2k+1$ it is enough to show a bijection between the vertex sets that maps non-adjacent vertices to non-adjacent vertices.

Let $X_1, \dots X_n$ be the vertices of $Q(n,k)$ in this order in the natural representation around the base cycle $C_n$. Let $f: V(Q(n,k)) \to V(K_{n/k})$ be defined by
$$X_u \mapsto uk \ ({\rm mod}\ n)$$
This is a one-to-one function since
$\gcd(n,k) = 1$.
Now look at $X_u \neq X_v$ arbitrary non-adjacent vertices in $Q(n,k)$. Let $\ell := |u-v|$ be their distance measured in rotations. If they are not adjacent, then one of them must be a right $j$-neighbor of the other for some $j \in \{1,\dots k-1\}$. Since all $j$-neighbors in $Q(n,k)$ have to be either $\ell$ or $\ell-1$ rotations apart, or they all have to be $\ell$ or $\ell+1$ rotations apart one of the equations $(k-x)\ell+x(\ell +1) = jn$ or $(k-x)\ell+x(\ell -1) = jn$ has an integral solution with $0<x<k$. (This is because if we consider the clockwise arc from each $z\in X_u$ to the $z'\in X_u$ for which this arc contains $j$ elements of $X_u$ including $z'$ but excluding $z$, then we cover $C_n$ exactly $j$ times.)
That means that $k\ell$ must belong to the same congruent class modulo $n$ as $x$ or $-x$, meaning that in the image the vertices $uk \ ({\rm mod}\ n)$ and $vk \ ({\rm mod}\ n)$ should be either less than $k$, or more than $n-k$ apart, i.e., they are indeed non-adjacent in $K_{n/k}$.
\qed

\begin{cor}\label{cor:qknk}
$$Q(n,k)\cong K_{n'/k'}$$
where $n'=\frac{n}{\gcd(n,k)}$ and $k'=\frac{k}{\gcd(n,k)}$.
\end{cor}

\proof
The statement is an immediate consequence of Propositions~\ref{prop:isotocirc} and \ref{prop:relpr}.
\qed

\medskip
\par\noindent
\begin{remark}
{\rm It is not hard to see that $K_{n/k}$ is always an induced subgraph of $\KG$ (see Proposition 6.19 in \cite{HN}), since putting the $n$ elements of the basic set on a cycle, the $k$-tuples of consecutive elements will be disjoint if and only if they are at least $k$ but at most $(n-k)$ rotations apart. This subgraph, however, may not appear in the Schrijver graph belonging to this defining cycle and it is not critical for the fractional chromatic number if $\gcd(n,k)>1$. Our results show that $\KG$ also contains $K_{n'/k'}$ as an induced subgraph for $n'=\frac{n}{\gcd(n,k)}, k'=\frac{k}{\gcd(n,k)}$, however, if we do not want to find this subgraph within a Schrijver graph then this can also be seen more directly. Put the elements of the basic $n$-set simply on $\ell:=\gcd(n,k)$ disjoint cycles of length $n'$ each and consider the $k$-subsets formed by the union of $k'$ consecutive elements on each of the $\ell$ cycles positioned similarly as above. More formally, if the elements of the $i^{\rm th}$ cycle are $x_{i,0},x_{i,1},\dots,x_{i,n'-1}$, then the $k$-subsets of the form
$$A_h=\{x_{i,h+j}\colon 1\le i\le \ell, 0\le j\le k'-1\},$$ (where $0\le h\le n'-1$ and addition is intended modulo $n'$)
induce a $K_{n'/k'}$ subgraph of $\KG$. Denoting the vertices of $K_{n'/k'}$ by $0,1,\dots,n'-1$ the isomorphism is given by the mapping $\varphi\colon h\mapsto A_h$.
\hfill$\Diamond$}
\end{remark}

\medskip
\par\noindent
\begin{remark}
{\rm The structure of $4$-chromatic Schrijver graphs ${\rm SG}(2k+2,k)$ is quite well-understood, cf. \cite{BBraun1, BBraun2, Li, STSzE80}. This structure is made of several levels, most of which are just cycles of length $2k+2$ with two exceptions at the two extremal levels. In one of those the cycle is extended to a complete bipartite graph while the other is extended to a M\"obius ladder in case $k$ is odd and is substituted by an odd cycle of length $k+1$ in case $k$ is even. This latter level is exactly the subgraph $Q(2k+2,k)$ induced by the well-spread $k$-subsets in this case.
\hfill$\Diamond$}
\end{remark}

\medskip
\par\noindent
\begin{prop}\label{prop:chiB}
  For all $n\ge 2k$ we have $$\chi(Q(n,k))=\left\lceil\frac{n}{k}\right\rceil.$$
\end{prop}

\proof
From Corollary~\ref{cor:qknk} and the properties of the circular complete graphs it follows that $$ \chi(Q(n,k))= \chi(K_{n'/k'}) =\left\lceil\frac{n'}{k'}\right\rceil=\left\lceil\frac{n}{k}\right\rceil,$$
where $n'=\frac{n}{\gcd(n,k)}, k'=\frac{k}{\gcd(n,k)}$.
\qed

\medskip
\par\noindent
One can also get to the same conclusion in a different way.
Following \cite{LPSV, STSzE80} we call an edge $XY\in E(\SG)$ interlacing if the elements of $X$ and $Y$ alternate on the defining cycle $C_n$. Let $I(n,k)$ denote the {\em interlacing graph} with parameters $n,k$ that is the subgraph of $\SG$ one obtains from $\SG$ after removing all edges that are not interlacing. It is proven in \cite{LPSV} and also follows from Theorem 6.33 in \cite{HN} that $$\chi(I(n,k))=\left\lceil\frac{n}{k}\right\rceil.$$

Note that it follows from Lemma~\ref{lem:wspalt} that all edges of $Q(n,k)$ are interlacing and this implies $Q(n,k)\subseteq I(n,k)$ which by the above mentioned result in \cite{LPSV, HN} gives that $\chi(Q(n,k))\le \chi(I(n,k))=\left\lceil\frac{n}{k}\right\rceil.$
On the other hand, by the first statement in Theorem~\ref{thm:main} we have $\chi(Q(n,k))\ge\chi_f(Q(n,k))=\frac{n}{k}$ and thus the integrality of $\chi(Q(n,k))$ gives the reverse inequality.

\medskip
\par\noindent
Note that Proposition~\ref{prop:chiB} gives a second proof for Corollary~\ref{cor:Brooks} as $Q(n,k)=\SG$ implies the equality of their chromatic number and $n-2k+2=\left\lceil\frac{n}{k}\right\rceil$ also implies that we must have $n=2k, n=2k+1$ or $k=1$.

\section{Critical edges}\label{sect:edges}

\medskip
\par\noindent
Here we are going to prove a strengthening of the second statement of Theorem~\ref{thm:main}.

\medskip
\par\noindent
\begin{defi}\label{defi:cedge}
For two adjacent vertices $X,Y\in V(Q(n,k))$ we call the edge $XY\in E(Q(n,k))$ a {\em cycle-edge} if the $k$-element set $X$ can be moved to $Y$ by rotating the defining cycle $C_n$ by just one element.
\end{defi}

\smallskip
\par\noindent
In other words, $XY\in E(Q(n,k))$ is a cycle-edge, if in the natural representation $X$ and $Y$ are next to each other along the base cycle.

\medskip
\par\noindent
\begin{thm}\label{thm:credge}
An edge of $Q(n,k)$ is critical for the fractional chromatic number if and only if it is a cycle-edge. More precisely, if $\gcd(n,k)=1$, $e\in E(Q(n,k))$ and $a,b$ are defined as in Theorem~\ref{thm:main} then
$$\chi_f(Q(n,k)\setminus\{e\})=
\left\{\begin{array}{lll} \frac{a}{b}&&\hbox{if $e$ is a cycle-edge}\\ \frac{n}{k}&&\hbox{otherwise.}\end{array}\right.$$
\end{thm}

\proof
As already mentioned it is a well-known fact that for every graph $G$ $\chi_f(G) \geq \frac{|V(G)|}{\alpha(G)}$. For $Q(n,k)$ we have $|V(Q(n,k))| = n$ and $\alpha(Q(n,k)) = k$, so if we want to decrease the fractional chromatic number by deleting an edge then we must find an edge whose absence allows to put an additional vertex into a maximum independent set.
This means that both endpoints of this edge should form a maximum independent set of the original graph with the rest of the vertices of the new maximum independent set.
In other words, we must find two independent sets of largest size both containing only one vertex not contained in the other.
\smallskip
\par\noindent
Recall that the maximum independents sets, that is those of size $k$ in $Q(n,k)$ are well-spread sets along the base cycle, so they are just rotations of each other.
\smallskip
\par\noindent
Now consider a maximum independent set $A_1$ in the natural representation.
We show that if we rotate $A_1$ along this cycle $a$ times (anticlockwise, say) then we get another maximum independent set $A_{a+1}$ for which $|A_1 \cap A_{a+1}| = k-1$ and thus $|A_1\cup A_{a+1}|=k+1$.
Recall from the proof of Lemma~\ref{lem:frszin} that if $A$ is a well-spread subset of size $k$ of an $n$-length cycle, where $a$ is the smallest positive solution of the diophantine equation $ak=bn-1$, then there is exactly one $a$-length arc of the cycle that contains only $b-1$ elements of $A$, all the other ones intersect $A$ in exactly $b$ elements. This means that all but one element of $A$ if moved by $a$ rotations along the cycle will hit another element of $A$. Choosing $A$ to be our maximum independent set $A_1$ this shows that $|A_1\cap A_{a+1}|\ge k-1$. This intersection cannot be larger, however, because that would mean $A_{a+1}=A_1$, but since $\gcd(n,k)=1$ we do have $n$ different rotations of $A_1$, so since $a<n$ we must have $A_{a+1}\neq A_1$.

\medskip
\par\noindent
Let $Y \in A_1  \setminus A_{a+1}$ and $X \in A_{a+1} \setminus A_1$ the two vertices in which $A_1$ and $A_{a+1}$ differ.
We know that $XY$ is an edge of $Q(n,k)$ (otherwise $A_1$ would have not been a maximum independent set). Now we show that deleting this edge
there exists a fractional coloring using $\frac{a}{b}$ total weight. We will modify the coloring used in the proof of Proposition~\ref{lem:frszin}.
Let $A_2, \dots, A_a$ be the independent sets we get by $1, \dots, a-1$ (anticlockwise) rotations of $A_1$.

\medskip
\par\noindent
The above argument implies that $X$ is the unique vertex on the base cycle that is covered only $b-1$ times by the independent sets $A_1,\dots,A_a$, all other vertices are covered by them $b$ times.
This means that if we put $\frac{1}{b}$ weight on the independent sets $A_1,\dots,A_a$ then every vertex but $X$ receives already a total weight $1$ (as the sum of the weights of independent sets containing it), while $X$ receives only $\frac{b-1}{b}$. But after deleting the edge $XY$ we can extend $A_1$ to $A_1':=A_1 \cup A_{a+1}=A_1\cup\{X\}$, and giving the $\frac{1}{b}$ weight of $A_1$ to $A_1'$ the total weight on $X$ will also be $1$. So this is a fractional coloring with total weight $\frac{a}{b}$ proving that $\chi_f(Q_{n,k}\setminus\{XY\})\le\frac{a}{b}$. Since the fractional chromatic number was $\frac{a}{b}$ even when deleting a vertex, this must hold with equality.

\medskip
\par\noindent
We still have to prove that $XY$ was necessarily a cycle edge and in fact, more generally, any edge the deletion of which from $Q(n,k)$ leads to a larger than $k$ independent set must be a cycle edge. For this it is enough to show that if $U$ and $W$ are two well-spread sets of size $k$ on an $n$ length cycle $C_n$ such that $|U\setminus W|=|W\setminus U|=1$ then the unique elements $u$ and $w$ of $U\setminus W$ and $W\setminus U$, respectively, should be neighboring elements of the cycle. We may also use that $k\le n/2$ so neither $U$ nor $W$ can contain consecutive elements of the cycle.
\smallskip\par\noindent
Assume for contradiction that the above is not true. Let $z\in U\cap W$ be the element which is closest on the cycle to $u$ on the arc which is the longer one between $u$ and $w$ and let $s$ be the neighbor of $u$ on the cycle in the other direction (i.e. towards $w$ on the arc where it is closer). By our assumptions $s\neq w$ and $s\notin U$. Now consider the arc $L$ starting with $z$ and ending with $u$ that gives $|L\cap U|=2$ by the choice of $z$. Let $L'$ be the first rotation of $L$ towards $w$. Thus $L'$ starts with the element following $z$ (towards $u$) and ends with $s\notin U$. Thus $L'\cap U=\{u\}$, in particular, $|L'\cap U|=1$ and since $w\notin L'$ it also implies $|L'\cap W|=0$. Since $W$ must be a rotated version of $U$ and $|L|=|L'|$ this contradicts the well-spreadness of $U$ completing the proof.
\qed

\medskip
\par\noindent
Similarly to the vertex-critical case, we can state a little stronger result.

\begin{prop}\label{prop:retractedge}
Assume $n\ge 2k$, $\gcd(n,k)=1$, let $a, b$ be the smallest positive integers for which $ak=bn-1$ and let $e$ be a cycle-edge of $Q(n,k)$. Then
$Q(n,k)\setminus\{e\}$
is homomorphically equivalent to $Q(a,b)$.
\end{prop}

\proof
Here the proof is essentially the same as that of Proposition~\ref{prop:retract}.
The homomorphism from $Q(a,b)$ to $Q(n,k)\setminus\{e\}$ we know to exist as before since $Q(a,b)$ is a subgraph of $Q(n,k)$ so again we only have to show that the reverse homomorphism exists.
We define the arc $B$ on the defining cycle $C_n$, the vertices $X_1, \dots, X_{n-1}$ of $Q(n,k)$, and $Y_1, \dots, Y_a$ of $Q(a,b)$ and the function $g(X_i) = Y_{i'}$, where $i \equiv i'$ (mod $a$) similarly as in the proof of Proposition~\ref{prop:retract}. Furthermore, let $X_n$ be the remaining vertex of $Q(n,k)$.
We still have for each $i \in [n-1]$ that $|X_i\cap B| = b$ and now we have $|X_n \cap B| = b-1$.
This means that in almost every $1$-element rotation from $X_i$ to $X_{i+1}$ the last vertex (clockwise endpoint) of $B$ will belong to $X_{i+1}$ if and only if the first vertex of $B$ belonged to $X_i$. There are only two  exceptions. The first one is the rotation from $X_{n-1}$ to $X_n$, where the first vertex of $B$ had to belong to $X_{n-1}$, but the last vertex cannot belong to $X_n$. The second is the rotation from $X_n$ to $X_1$ where the number of vertices belonging to $B$ has to increase. So $\{g(v) | v\in X_n \cap B\}$ will not be the same as $g(X_n)$.
In the case of $X_1$ and $X_n$ we will have that $X_1 \cap X_n = \emptyset$, but $g(X_1) \cap g(X_n) \neq \emptyset$. So $\{X_1,X_n\}\in E(Q(n,k))$ while $\{g(X_1),g(X_n)\}\notin E(Q(a,b))$. But $\{X_1,X_n\}$ is a cycle-edge of $Q(n,k)$ so (by symmetry) we may choose it to be the deleted edge $e$ making $g$ a homomorphism from $Q(n,k)\setminus\{e\}$ to $Q(a,b)$.
\qed

\medskip
\par\noindent

The circular chromatic number $\chi_c$ is a graph parameter that can be defined via the existence of graph homomorphisms to circular complete graphs (cf. e.g. \cite{HN} for more details). We know that $\chi_c(K_{n/k}) = \frac{n}{k}$ and it has some interest to see how this value may change if we remove an edge from $K_{n/k}$. By the previous results this can be answered as follows.

\begin{cor}
If $\gcd(n,k)=1$, $e :=\{i,j \}\in E(K_{n/k})$ with $j>i$ and $a,b$ are defined as in Theorem~\ref{thm:main} then
$$\chi_c(K_{n/k})\setminus\{e\})=
\left\{\begin{array}{lll} \frac{a}{b}&&\hbox{if $j-i = k\ {\rm or}\ n-k$}\\ \frac{n}{k}&&\hbox{otherwise.}\end{array}\right.$$
\end{cor}
\proof
We know from Proposition~\ref{prop:isotocirc} that $Q(m,r)$ is isomorphic with the circular complete graph $K_{m/r}$ when $\gcd(m,r)=1$ thus we have both $Q(n,k)\cong K_{n/k}$ and $Q(a,b)\cong K_(a/b)$ observing that $\gcd(n,k)=1$ and the definition of $a,b$ also implies $\gcd(a,b)=1$.
\smallskip
\par\noindent
With this in mind the statement is just a consequence of Theorem~\ref{thm:credge}, Proposition~\ref{prop:retractedge}, the fact, that for every graph $G$ we have $\chi_f(G) \leq \chi_c(G)$ and the observation that the isomorphism given in the proof of Proposition~\ref{prop:isotocirc} maps the endvertices of cycle edges of $Q(n,k)$ to vertices having (circular) distance $k$ in $K_{n/k}$.
\qed

\end{document}